\documentclass[11pt]{article}
\usepackage{amsmath,amsthm,amsfonts,a4,bezier}

\usepackage{amssymb}
\usepackage{hyperref}
\usepackage{amscd}
\usepackage{graphpap}
\usepackage[pdftex]{color,graphicx,xcolor}
\usepackage{color,graphicx,enumerate}
\usepackage{setspace}
\usepackage{comment}
\usepackage{mathtools}
\usepackage{relsize}
\usepackage{indentfirst}
\usepackage{verbatim,bbm}
\usepackage{hyperref}
\usepackage{extarrows}   
\usepackage[pdftex]{color}
\usepackage{relsize}
\usepackage{calrsfs}
\usepackage{centernot}
\usepackage{tikz}

\DeclareGraphicsExtensions{.jpg,.pdf,.eps,.png}
\newtheorem{theorem}{Theorem}

\newtheorem{proposition}{Proposition}
\newtheorem{lemma}{Lemma}

\newtheorem{remark}{Remark}

\renewcommand{\leq}{\leqslant}
\renewcommand{\geq}{\geqslant}

\newcommand{\E}{\mathcal{E}}

\newcommand{\Z}{{\mathbb{Z}}}
\newcommand{\LL}{{\mathbb{L}}}
\newcommand{\N}{{\mathbb{N}}}

\renewcommand{\P}{{\mathbb{P}}}

\newcommand{\V}{{\mathcal {V}}}
\newcommand{\F}{{\mathcal {F}}}

\newcommand{\G}{{\mathcal{G}}}
\begin{document}
\thispagestyle{empty}
\baselineskip=28pt
\vskip 5mm
\begin{center} {\LARGE{\bf Percolation of words on the hypercubic lattice with one-dimensional long-range interactions}}
\end{center}

\baselineskip=15pt
\vskip 10mm

\begin{center}\large
Pablo A. Gomes\footnote{Instituto de Matem\'atica e  Estat\'{\i}stica, Universidade de S\~ao Paulo, Brazil, {pagomes@usp.br}}, Otávio Lima\footnote{Departamento de Estat\'{\i}stica, Universidade Federal de Minas Gerais, Brazil, {otaviooasl@ufmg.br}} and  Roger W. C. Silva 
 \footnote{Departamento de Estat\'{\i}stica,
 Universidade Federal de Minas Gerais, Brazil, \href{mailto:rogerwcs@ufmg.br} {rogerwcs@est.ufmg.br;}}
\\ 


\end{center}
\vspace{1cm}
\begin{abstract} 
We investigate the problem of percolation of words in a random environment. To each vertex, we independently assign a letter $0$ or $1$ according to Bernoulli r.v.'s with parameter $p$. The environment is the resulting graph obtained from an independent long-range bond percolation configuration on $\Z^{d-1} \times \Z$, $d\geq 3$, where each edge parallel to $\Z^{d-1}$ has length one and is open with probability $\epsilon$, while edges of length $n$ parallel to $\Z$ are open with probability $p_n$. We prove that if the sum of $p_n$ diverges, then for any $\epsilon$ and $p$, there is a $K$ such that all words are seen from the origin with probability close to $1$, even if all connections with length larger than $K$ are suppressed.
\medskip


\end{abstract}

\onehalfspacing
\section{Introduction}

The problem of percolation of words was introduced in a seminal paper by I. Benjamini and H. Kesten (see \cite{BK}) and is formulated as follows: let $\G=(\V,\E)$ be a graph with vertex set $\V$, which is assumed to be countably infinite. To each site $v\in\V$ we assign independently a random variable $X(v)$, which takes the value 1 or 0 with probability $p$ or $1-p$, respectively. This induces the probability space $(\Omega_S,\F,P_p)$, where $\Omega_S=\{0,1\}^{\V}$, $\F$ is the sigma-field generated by the cylinder sets of $\Omega_S$ and $P_p$ is the product over $\V$ of the measures assigning masses $p$ and $1-p$ to the points 1 and 0, respectively. A typical element of $\Omega_S$ is denoted by $\omega_S$ and we write $X(v,\omega_S)$ for the state of vertex $v$ in the configuration $\omega_S$. A \textit{path} on the graph $\G$ is a sequence $(v_1,v_2,\dots)$, $v_i\in\V$, $v_i\neq v_j$ for $i\neq j$, and such that $\langle v_i,v_j\rangle$ is an element of $\E$.

Let $$\Xi=\{0,1\}^{\mathbb{N}}.$$ 
An element of $\Xi$ is denoted by $\xi=\{\xi_n\}_{n\in\mathbb{N}}$ and is called a \textit{word}.  We say $\xi$ is \textit{seen from the vertex} $v\in\V$ in the configuration $\omega_S$ if there is a path $(v,v_1,v_2,\dots)$ on $\G$ such that $X(v_i,\omega_S)=\xi_i$, $i\geq 1$. Note that the state of $v$ plays no role here.  

Write
$$W_v=W_v(\omega_S)=\{\xi\in \Xi: \xi \mbox{ is seem from $v$ in }\omega_S\},$$
and
$$W_{\infty}=W_{\infty}(\omega_S)=\bigcup_{v\in\V}W_v(\omega_S),$$ the collection of words seen from some vertex in $\G$. Clearly the largest these sets can be is $\Xi$ and, to the best of our knowledge, the first one to ask when $W_v$ can in fact equal all of $\Xi$ was F. Dekking (see \cite{D}).

The model we just described is a natural generalization of the usual percolation framework (classical Bernoulli percolation occurs when the word $(1,1,\dots)$ is seen from some vertex $v\in\V$), and consequently, it is a source of interesting mathematical questions. Notably, the model received much attention in the last three decades and a number of papers were written about the subject. We do not intend to give a complete review of the results obtained so far, but rather a brief description of some of them. In what follows, $\LL^d$ stands for the usual hypercubic lattice with nearest neighbors.    

In \cite{BK}, the authors study the problem of seeing all words on the graph $\LL^d$ and show that the event $\{W_{\infty}=\Xi\}$ occurs almost surely when $p=1/2$ and $d\geq 10$. Also, the event $\{W_{v}=\Xi\mbox{ for some $v$}\}$ occurs almost surely when  $p=1/2$ and $d\geq 40$. The authors make it clear that those bounds on $d$ are not sharp, albeit their methods did not allow for an improvement.  Three years later, H. Kesten, V. Sidoravicius, and Y. Zhang \cite{KSZ} showed that, with probability 1, \textit{almost all} words can be seen in critical site percolation on the triangular lattice, answering Open Problem 1 in \cite{BK}. Here, almost all refers to the product measure $\nu_{\lambda}\coloneqq ((1-\lambda)\delta_0+\lambda\delta_1))^{\otimes\N}$ on the set of words $\Xi$, $0<\lambda<1$. In a second paper  (see \cite{KSZ2}), the same authors investigate the problem of seeing all words in site percolation on the closed-packed graph of the square lattice, which is obtained by adding both diagonal edges to each face of $\LL^2$.  They show that for every fixed $p\in(1-p_c(\LL^2),p_c(\LL^2))$, the event $\{W_{v}=\Xi\mbox{ for some $v$}\}$ has probability 1. Here $p_c(\G)$ stands for the critical threshold of Bernoulli site percolation on a graph $\G$.

It is not hard to see that $\nu_{1/2}-$almost all words can be seen in any graph $\G$ with $p_c(\G)<1/2$ when $p=1/2$. This follows by Wierman`s coupling (see \cite{W} and the discussion following Open Problem 1 in \cite{BK}). In particular, since $p_c(\LL^3)<1/2$ (see \cite{CR}), $\nu_{1/2}-$almost all words can be seen in $\LL^3$ when $p=1/2$. Open Problem 2 in \cite{BK} asks if the stronger claim that all words can be seen in $\LL^3$ when $p=1/2$ is true. In a recent paper (see \cite{NTT}), the authors give a positive answer to this question. In fact, they prove the stronger statement that all words can be seen from some vertex in a sufficiently thick slab of $\LL^d$, for all $d\geq 3$ and $p\in(p_c(\LL^d),1-p_c(\LL^d))$.

We observe that in general, the problem of seeing all words is harder than the problem of seeing almost all words. See Section 7 of \cite{BK} for an example where one sees almost all but not all words.

\subsection{The long-range setting}
In this work, we study the problem of percolation of words on the hypercubic lattice with long-range connections in one direction. However, unlike all previous works related to this problem, we do not have all edges of the graph present in order to find a path that coincides with a given word. Instead, we try to see words in the resulting random graph from a long-range bond percolation configuration. More precisely, let $(\mathbf{e}_1,\dots,\mathbf{e}_d)$ be the canonical basis of $\mathbb{R}^d$, and define the sets 
$$\E_V=\{\langle u,u+n\mathbf{e}_d\rangle: u\in \Z^d, n\in\N\},$$
$$\E_H=\{\langle u,u+\mathbf{e}_i\rangle: u\in \Z^{d},  i=1,\dots,d-1\},$$
$$\E=\E_V\cup \E_H.$$

We consider the following oriented bond percolation model on the graph $\mathbb{G}^d=(\Z^d,\E_V\cup \E_H)$. Given a sequence $\{p_n\}_{n\in\N}$, $p_n\in[0,1]$, and $\epsilon>0$, consider a family $\{Y(e)\}_{e\in\E}$ of independent Bernoulli random variables with distribution $\{\mu_e\}_{e\in\E}$ such that
\begin{equation*}\mu_{\langle u, v \rangle}(Y(\langle u, v \rangle)=1)=\left\{\begin{array}{ll}
p_{||u-v||}&\mbox{if}\quad \langle u,v\rangle\in \E_V,\\
 \epsilon&\mbox{if}\quad \langle u,v\rangle\in\E_H,
\end{array}\right.
\end{equation*}
where $||u-v||$ denotes the graph distance between vertices $u$ and $v$. This induces the probability space $(\Omega_B,\mathcal{H},P_{\epsilon})$, where $\Omega_B=\{0,1\}^{\E}$, $\mathcal{H}$ is the usual sigma-field, and \begin{equation}\label{prod_meas}P_{\epsilon}=\prod_{\langle u,v \rangle \in\E}\mu_{\langle u,v \rangle}.
\end{equation} A typical element of $\Omega_B$ is denoted by $\omega_B$, and we write $Y(e, \omega_B)$ for the state of edge $e$ in the configuration $\omega_B$. Our main assumption on the sequence $\{p_n\}_{n\in\N}$ will be 
\begin{equation}\label{diverge}
\sum_{n=1}^{\infty}p_n=\infty.
\end{equation}

It is not hard to see that assumption \eqref{diverge} implies that all words can be seen from the origin almost surely. A natural question is the following: is the infinitude of long-range connections crucial for the occurrence of such event? 
With this question in mind, define the $K$-truncated sequence
\begin{equation*}
p_n^K=\left\{\begin{array}{ll}
p_n&\mbox{if}\quad n\leq K,\\
 0&\mbox{if}\quad n>K,
\end{array}\right.
\end{equation*}
$K\in\N$, and denote by $P_{\epsilon}^K$ the $K$-truncated product measure defined as in \eqref{prod_meas} with $p_n^K$ instead of $p_n$.

As before, let $\{X(v)\}_{v\in\Z^d}$ be a collection of Bernoulli $i.i.d.$ random variables with parameter $p$ and corresponding product measure $P_p$. We say a word $(\xi_1,\xi_2,\dots)$ \textit{is seen from the vertex $v\in\Z^d$} (in the $K$-truncated model) in the configuration $\omega_B\times\omega_S$ if there is a sequence $(v= v_0,v_1,v_2,\dots)$, $v_i\in\Z^d$, such that
$v_i\neq v_j$, $e_i=\langle v_{i-1},v_i\rangle \in \E$, $||v_{i-1},v_i||\leq K$, $X(v_i,\omega_S)=\xi_i$, $Y(e_i, \omega_B)=1$, for all $i,j$. If this occurs, we will also say that $\xi$ percolates in the $K$-truncated model.

It is clear that the above construction has two sources of randomness. First, a long-range bond configuration $\omega_B$ is sampled according to $P_{\epsilon}^K$, and then it is checked if $\xi$ is seen in $\omega_B$ according to an appropriate conditional measure $P_p^{\omega_B}$. Formally, let $\P_{p,\epsilon}^K$ denote the probability measure on the measurable space $(\Omega_B\times\Omega_S,\sigma(\mathcal{H}\times\mathcal{F}))$, such that, for all measurable rectangles of the form $R_1\times R_2$,
\begin{equation*}
\P_{p,\epsilon}^K(R_1\times R_2)=\displaystyle\int_{R_1}P_p^{\omega_B}(R_2)dP_{\epsilon}^K(\omega_B),\,\,\,\,\,\,\,R_1\in\mathcal{H}, R_2\in\mathcal{F},
\end{equation*}
where $P_p^{\omega_B}$ is the (quenched) law of independent site percolation with parameter $p$, conditional on the bond percolation configuration $\omega_B$.

There is a growing literature on long-range percolation models. We refer the reader for the introductory sections of \cite{ELV} and \cite{FL} for a comprehensive list of papers about the subject. Let us mention in more detail the results in \cite{ELV} and \cite{LSS}. In the former, the authors consider oriented percolation on $\LL^d$ with long-range connections in all directions. They show that, under the hypothesis in \eqref{diverge}, the word $(1,1,1,\dots)$ is seen from the origin with positive probability (see Theorem 1 in \cite{ELV}). In the same paper, the authors consider non-oriented percolation on $\LL^d$, with nearest neighbor edges (which are open with probability $\epsilon$) in the vertical direction and long-range edges in the remaining $d-1$ directions.  Again, they show the word $(1,1,1,\dots)$ is seem from the origin with positive probability if one assumes the condition in \eqref{diverge} (see Theorem 2 in \cite{ELV}). In \cite{LSS}, the problem of percolation of words is addressed, albeit with a more restricted hypothesis. The authors prove that, given any $p\in(0,1)$, there exists $K(p)\in\N$ such that all words can be seen from the origin with positive probability on $\LL^d$, $d\geq 2$, with long-range connections in all directions, and $p_n=1$ for all $n$.  

In the remainder of this paper, we investigate the problem of percolation of words in the long-range setting described at the beginning of this section. We will show that, under the hypothesis in \eqref{diverge}, the infinitude of long-range connections is not crucial for the occurrence of the event that all words are seen from the origin. Indeed, we will show that, given any $p\in(0,1)$ and any $\epsilon>0$, there is a $K(\{p_n\},p,\epsilon)$ large enough such that all words can be seen from the origin with probability close to 1. This result (see Theorem \ref{words} below) generalizes the main theorem in \cite{ELV}.

\section{Main result}

Our main contribution is Theorem \ref{words} below. We will prove it assuming an auxiliary result (see Proposition \ref{aux}) that will be shown in Section \ref{aux_proof}. In what follows, $\mathbb{G}^d_+$ denotes the subgraph of $\mathbb{G}^d$ with vertex set $\mathbb{Z}^d_+$ and edge set consisting of the edges in $\E_V\cup\E_H$ with both endpoints in $\Z^d_+$.

\subsection{Statement}

\begin{theorem}\label{words} Consider a long-range oriented percolation process on $\mathbb{G}^d_+$, $d\geq 3$, and assume $\{p_n\}_{n\in\N}$ satisfies \eqref{diverge}. Then, for all $p\in(0,1)$, $\epsilon>0$, and $\alpha>0$, there exists $K=K(\{p_n\},p,\epsilon,\alpha)\in\N$ large enough such that $$\P_{p,\epsilon}^K(W_0=\Xi)>1-\alpha.$$
\end{theorem}

\begin{remark}
By translation invariance, Theorem \ref{words} implies that $\P_{p,\epsilon}^K(W_v=\Xi\mbox{ for some $v$})=1.$ Therefore, we conclude that the environment $\omega_B$ is such that $$P_p^{\omega_B}(W_v=\Xi\mbox{ for some v})=1,\,\,\,P_{\epsilon}^K-a.s..$$  
\end{remark}

Note that it suffices to prove the theorem when $d=3$. Part of the proof consists in a dynamical coupling between the long-range model and an independent nearest neighbor oriented percolation process on the quadrant $\LL^2_+$ with parameter $\gamma\in(0,1)$, which will be large, in general. We describe the coupling in detail in the next section.

\subsection{Tools and proof of Theorem \ref{words}}\label{main_proof}

\subsubsection{The coupling}\label{coupling}

For a given set $A\subset \Z^2_+$, write $\partial_e A$ for the external vertex boundary of $A$, that is,
$$\partial_e A=\{u\in A^c\colon  \exists v\in {A}\mbox{ such that } ||u-v||=1\}.$$
Fix a word $\xi\in\Xi$. We will construct inductively a sequence $\{A_n,B_n\}_{n\geq 0}$ of ordered pairs of subsets of $\Z_+^2$ and a function 
$$\psi:A_n\longrightarrow \Z_+.$$ 
Consider an arbitrary, but fixed, ordering of the vertices of $\Z^2_+$. Let $\mathcal{O}=(0,0)$ be the origin of $\LL^2_+$ and set $A_0=\{\mathcal{O}\}$, $B_0=\emptyset$, and $\psi(\mathcal{O})=0$. Assume $\{A_n,B_n\}$ has been constructed for some $n\in\N$ and that $\psi(x)$ is known for all $x\in A_n$. If $\partial_e A_n\cap B_n^c=\emptyset$, then stop the exploration and set $(A_{\ell},B_{\ell})=(A_n,B_n)$, for all $\ell\geq n$. If $\partial_e A_n\cap B_n^c\neq\emptyset$, let $x_n$ be the earliest vertex in the fixed ordering in $\partial_e A_n\cap B_n^c$ and define $y_n$ as the vertex in $A_n$ such that $x_n=y_n+(1,0)$ or $x_n=y_n+(0,1)$. Assume first that $x_n=y_n+(0,1)$. For a fixed $N\in\N$, we say $x_n$ is \textit{black} if, for some $i\in\{1,\dots, N\}$, all the following conditions hold:
\begin{itemize}
\item $X(y_n,\psi(y_n)+i)=\xi_{2||y_n||+1}$, $(y_n,\psi(y_n)+i)\in \Z_+^2\times \Z_+,$
\item $X(x_n,\psi(y_n)+i)=\xi_{2||y_n||+2}$, $(x_n,\psi(y_n)+i)\in \Z_+^2\times \Z_+$,
\item $\langle (y_n,\psi(y_n)),(y_n,\psi(y_n)+i)\rangle \in \E_V$, and $\langle (y_n,\psi(y_n)+i),(x_n,\psi(y_n)+i)\rangle \in \E_H$ are open.
\end{itemize}
Write $\psi(x_n)=\psi(y_n)+i$ and
define
\begin{equation*}(A_{n+1},B_{n+1})=\left\{\begin{array}{ll}
(A_n\cup\{x_n\},B_n)&\mbox{if}\quad x_n\mbox{ is black},\\
 (A_n,B_n\cup\{x_n\})&\mbox{otherwise}.
\end{array}\right.
\end{equation*}
If $x_n=y_n+(0,1)$, we proceed analogously with $i\in\{N+1,\dots,N+M\}$. Such distinction is useful to avoid dependence issues. If there are two $y_n's$ available we choose the one of the first kind. Writing $$A(\xi)=\bigcup_{n\in\N}A_n,$$ it is clear that if $|A(\xi)|=\infty$, then $\xi$ is seen in the $(N+M)$-truncated model. Note that
\begin{equation*}\P_{p,\epsilon}^{N+M}(x_n\mbox{ is not black})\leq\left\{\begin{array}{ll}
\prod_{i=1}^N[1-\epsilon p_i(\min\{p,1-p\})^2]&\mbox{if}\quad x_n=y_n+(1,0),\\
 \prod_{i=N+1}^{N+M}[1-\epsilon p_i(\min\{p,1-p\})^2]&\mbox{if}\quad x_n=y_n+(0,1).
\end{array}\right.
\end{equation*}
In any case, assumption \eqref{diverge} implies
\begin{equation}\label{black}
\P_{p,\epsilon}^{N+M}(x_n\mbox{ is black})\longrightarrow 1,\mbox{ when $N,M\rightarrow \infty$},
\end{equation}
for all $n\in\N$. 

A key feature of the coupling above is that for any $p$, $\epsilon$, and $0<\gamma<1$, we can choose $N=N(\{p_n\},p,\epsilon,\gamma)$ and $M=M(\{p_n\},p,\epsilon,\gamma)$ sufficiently large so that $\xi$-percolation on the $(N+M)$-truncated model (or the corresponding process of black points) dominates an independent oriented percolation process with density $\gamma$. This is guaranteed by Lemma 1 in \cite{GM} and shall be useful later, when control over a highly supercritical oriented percolation process on $\LL^2_+$ will be needed. 

\subsubsection{Auxiliary result}
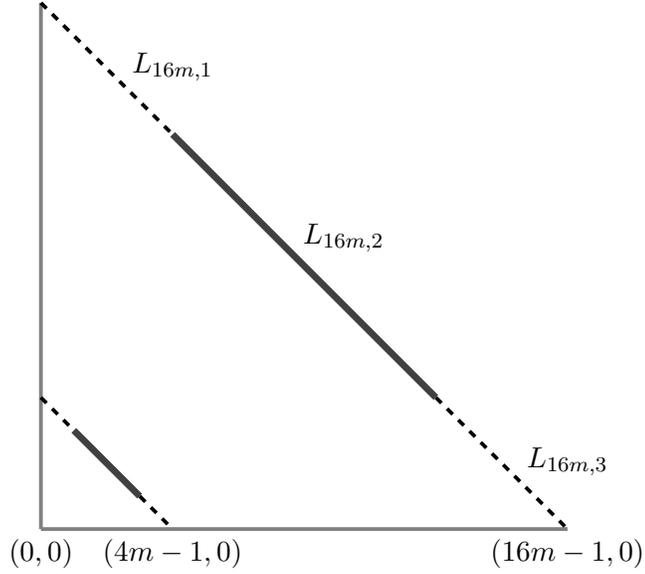
\begin{figure}[ht]
\centering
\begin{tikzpicture}[scale=0.035]

\draw[line width=0.05cm,gray] (0,0) -- (0,200);
\draw[line width=0.05cm,gray] (0,0) -- (200,0);
\draw[line width=0.05cm,dashed] (0,200) -- (50,150);
\draw[line width=0.05cm,dashed] (150,50) -- (200,0);
\draw[line width=0.09cm,darkgray] (50,150) -- (150,50);

\draw[line width=0.05cm,dashed] (0,50) -- (12.5,37.5);
\draw[line width=0.09cm,darkgray] (12.5,37.5) -- (37.5,12.5);
\draw[line width=0.05cm,dashed] (37.5,12.5) -- (50,0);

\draw (0,0)node[below] {$(0,0)$};
\draw (200,0)node[below] {$(16m-1,0)$};
\draw (50,0)node[below] {$(4m-1,0)$};
\draw (115,120)node[below] {\normalsize $L_{16m,2}$};
\draw (50,185)node[below] {\normalsize $L_{16m,1}$};
\draw (200,35)node[below] {\normalsize $L_{16m,3}$};


\end{tikzpicture}
\caption{Graphical representation of the sets $L_{16m,1}$, $L_{16m,2}$, and $L_{16m,3}$\label{rep_1}.}
\label{Lambda}
\end{figure}

Before we prove Theorem \ref{words} we need to introduce further notation.
For each $n\in\N$, let $$\Xi_n=\{0,1\}^{2(n-1)}$$ be the set of finite words $(\xi_1,\xi_2,\dots,\xi_{2(n-1)})$ of length $2(n-1)$. Define the function
$$\sigma_n:\left(\bigcup_{m\geq n}\Xi_m\right)\cup \Xi\longrightarrow \Xi_n,$$
where $\sigma_n(\xi_1,\dots,\xi_{2m-2})=(\xi_1,\dots,\xi_{2(n-1)})$, and $\sigma_n(\xi_1,\xi_2,\dots)=(\xi_1,\dots,\xi_{2(n-1)})$.
For any $m=4\ell$, $\ell\in\N$, write
\begin{equation*}
L_{m}=\{v\in\Z^2_+:||v||=m-1\},
\end{equation*}
and consider the line segment passing by the points in $L_{m}$. We partition $L_{m}$ in the sets $L_{m,r}$, $r=1,2,3$, with
\begin{align*}
    L_{m,1} &=\{(v_1,v_2)\in L_{m}: 0\leq v_1 < m/4\}, \\
    L_{m,2} &=\{(v_1,v_2)\in L_{m}: m/4\leq v_1 < 3m/4\},\\
    L_{m,3} &=\{(v_1,v_2)\in L_{m}: 3m/4 \leq v_1 < m\}.
\end{align*}
See Figure \ref{rep_1} for a graphical representation of the sets above.

For each $\eta\in \Xi_{4m}$, pick any $\xi\in \Xi$ such that $\sigma_{4m}(\xi)=\eta$ and define $$\Gamma_m(\eta)=A(\xi)\cap L_{4m,2}.$$
The reason to choose finite words of length $2(n-1)$ is because at each step of the coupling construction we see exactly two digits of the word. Observe that $\Gamma_m(\eta)$ does not depend on the choice of $\xi$, since $A(\xi)\cap L_{4m,2}$ depends only on the first $8m-2$ digits of $\xi$. Define 
$$\mathcal{B}_m(\eta)=\{|\Gamma_m(\eta)|\geq m\},$$
and write 
$$\mathcal{D}_{m}=\bigcap_{\eta\in \Xi_{4m}}\mathcal{B}_{m}(\eta).$$

The next proposition is inspired by \cite{NTT}. 
Intuitively speaking, it says that conditionally on the event
that at least a fraction of 1/2 of the vertices of $L_{4m,2}$ are black (for a fixed word $\eta$), then it is very unlikely (with a suitable choice of $N$ and $M$) that the fraction of black points on $L_{16m,2}$ gets below 1/2. Such a result will be important to balance the probability of seeing a word with the entropy created by the
words.
\begin{proposition}\label{aux} Let $p\in(0,1)$, $\epsilon>0$, and $\{p_n\}_{n\in\N}$ such that \eqref{diverge} holds. For any $a>0$, there are $N=N(\{p_n\},a,p,\epsilon)$ and $M=M(\{p_n\},a,p,\epsilon)$ sufficiently large such that $$\P_{p,\epsilon}^{N+M}(\mathcal{B}_{4m}^c(\eta)\cap\mathcal{B}_{m}(\sigma_{4m}(\eta)))\leq a^{m},$$ for all $m\geq 1$, and for all $\eta\in \Xi_{16m}$. 
\end{proposition}

The proof of Proposition \ref{aux} is postponed to the next section. Let us prove Theorem \ref{words}.

\hspace{-0.6cm}\textit{Proof of Theorem \ref{words}.} The first thing to note is that 
\begin{equation*} \{W_{0}=\Xi\}\supset \bigcap_{m\in \Lambda}^{\infty}\mathcal{D}_{m},
\end{equation*}
where $\Lambda=\{4^j:j\in\N\}$.
Writing
$$\left(\bigcap_{m\in\Lambda}\mathcal{D}_m\right)^c=\bigcup_{\substack{m\in\Lambda\\m\neq 1}}\left(\mathcal{D}_{4m}^c\cap\mathcal{D}_m\right)\cup \mathcal{D}_1^c,$$
and observing that \eqref{black} implies
$$\lim_{N,M\rightarrow \infty}\P_{p,\epsilon}^{N+M}(\mathcal{D}_1^c)=0,$$ it suffices to show that
$$\lim_{N,M\rightarrow \infty}\sum_{\substack{m\in\Lambda\\m\neq 1}}\P_{p,\epsilon}^{N+M}(\mathcal{D}_{4m}^c\cap\mathcal{D}_m)\leq \lim_{N,M\rightarrow \infty}\sum_{m\in\N}\P_{p,\epsilon}^{N+M}(\mathcal{D}_{4m}^c\cap\mathcal{D}_m)=0.$$

By definition of $\mathcal{D}_m$, we have
\begin{equation*}
\mathcal{D}_{4m}^c\cap\mathcal{D}_m=\bigcup_{\eta\in \Xi_{16m}}\left\{\mathcal{B}_{4m}^c(\eta)\cap\mathcal{B}_{m}(\sigma_{4m}(\eta)\right)\}.
\end{equation*}
Taking $N$ and $M$ as in Proposition \ref{aux} we obtain
$$\sum_{m\in\N}\P_{p,\epsilon}^{N+M}\left(\mathcal{D}_{4m}^c\cap\mathcal{D}_m\right)\leq \sum_{m\in\N}|\Xi_{16m}|a^m \leq \sum_{m\in\N}(2^{32}a)^m.$$
Since $a$ is arbitrary, the result follows with $K=N+M$.
\qed

\subsection{Proof of Proposition \ref{aux}}\label{aux_proof}

The idea of the proof is to use the dynamical coupling described in Section \ref{coupling} to compare the truncated long-range percolation model on $\mathbb{G}^3_+$ with a highly supercritical oriented percolation process on $\LL^2_+$ (this type of comparison is reminiscent from the '80s and '90s, see for instance \cite{CCN} and \cite{GS}). Such coupling combined with a high density ``propagation''  of the oriented percolation cluster on $\LL^2_+$ will give the desired result. Let us get into the details.

Let $\gamma\in(0,1)$ and consider an independent nearest neighbor oriented site percolation process on the quadrant $\LL^2_+$ with parameter $\gamma$. That is, to each site $v\in\Z^2_+$ we assign independently a Bernoulli random variable $\omega_v$ with parameter $\gamma$. If $w_v=1$ ($w_v=0$) we say $v$ is \textit{occupied} (\textit{vacant}). Let $\P_{\gamma}$ denote the corresponding product law.

We denote by $\{x\rightarrow y\}$ the event where $x\in\Z^2_+$ is connected to $y\in\Z^2_+$ by an oriented path of occupied vertices. Also, for any $E,F\subset \Z^2_+$, write $\{E\rightarrow F\}=\{\exists\,x\in E,y\in F: x\rightarrow y \}$. 

For each $m\in\N$ and $S\subset L_{4m}$, let $$M_S=|\{x\in L_{16m,2}:S\rightarrow x\}|.$$
The proof of Proposition \ref{aux} follows from the following two lemmas whose proofs are given in the next section.
\begin{lemma}\label{orien_perc} For every $a>0$, there exist $\gamma(a)\in(0,1)$ such that $$\P_{\gamma}\left(M_S<4m\right)\leq a^{m},$$
$\forall m\in\N$, $\forall S\subset L_{4m,2}$ such that $|S|\geq m$. 
\end{lemma}

\begin{lemma}\label{compar} Let $\gamma\in[0,1)$, $p\in(0,1)$, $\epsilon>0$, and $\{p_n\}_{n\in\N}$ such that \eqref{diverge} holds. There are $N=N(\{p_n\},p,\epsilon,\gamma)\in\N$ and $M=M(\{p_n\},p,\epsilon,\gamma)\in\N$ large enough such that 
$$\P_{p,\epsilon}^{N+M}\left(|\Gamma_{4m}(\eta)|<n\Bigm|\Gamma_m\left(\sigma_{4m}(\eta)\right)=S\right)\leq \P_{\gamma}\left(M_S<n\right),$$
$\forall n\in\N$, $\forall \eta\in \Xi_{16m}$, $\forall m\in\N$, and $\forall S\subset L_{4m,2}$.
\end{lemma}

To see why Proposition \ref{aux} follows from the lemmas above, simply note that for every $a>0$, there are $\gamma=\gamma(a)\in(0,1)$, $N=N(\{p_n\},p,\epsilon,\gamma)\in\N$ and $M=M(\{p_n\},p,\epsilon,\gamma)\in\N$ large enough such that
 \begin{align*}\P_{p,\epsilon}^{N+M}(\mathcal{B}_{4m}^c(\eta)\cap\mathcal{B}_{m}(\sigma_{4m}(\eta))) &= \sum_{\substack{S\subset L_{4m,2}\\|S|\geq m}}\P_{p,\epsilon}^{N+M}\left(|\Gamma_{4m}(\eta)| < 4m\cap\Gamma_m(\sigma_{4m}(\eta))=S\right)\\
 &\leq \sum_{\substack{S\subset L_{4m,2}\\|S|\geq m}} \P_{\gamma}(M_S<4m) \P_{p,\epsilon}^{N+M}(\Gamma_m(\sigma_{4m}(\eta))=S)\leq a^m.
 \end{align*}

\section{Proofs of Lemmas}

In this section, we prove Lemma \ref{orien_perc} and Lemma \ref{compar}. We shall need a result on the stochastic domination of oriented site percolation measure over product measure, for suitable choices of the density of open sites in both models.  Such result is due to Liggett and Steif (see Theorem~1.1 and the remark following Theorem 2.1 in \cite{LS}) and is originally formulated for the contact process. Before we dive into the proofs of the lemmas, let us state the theorem in the context of oriented site percolation. For that purpose, we need to introduce further notation.

Consider the independent oriented site percolation on $\LL^2$ with parameter $\gamma$. For each $t \in \Z$, let $V_t = \{ (x,y) \in \Z^2 : x+y = t\}$. Given $A \subset \Z$, let $A =  \{ (x,y) \in V_0 : x \in A\}$. For any $A' \subset \Z$ and $t\in\Z^+$, define
$$\xi_t^A=\{x\in\Z : A' \rightarrow (x,t-x) \in V_t\}.$$
Then, for any $A\subset\Z$, the process $\{\xi_t^{A}\}_{t\geq 0}$ is a time-homogeneous Markov chain, taking values in the subsets of $\Z$, with initial configuration $\xi_0^{A}=A$. Observe that
$$P_t(A,\cdot)\coloneqq \P_{\gamma}(\xi_{t}^{A}\in\cdot)$$ 
defines a transition probability on $\{0,1\}^{\Z}$. We call a Markov chain in $\{0,1\}^{\Z}$ with such transition probability, and arbitrary initial law, an \textit{oriented percolation process}. If one denotes by $-\infty \rightarrow (x,y)$ the event that the set $\{(u,v)\in\Z^2 : (u,v)\rightarrow (x,y)\}$ is infinite, and 
$$\xi_t^{\{-\infty\}}\coloneqq\{x\in\Z:-\infty\rightarrow (x,t-x)\}, ~t\in\Z,$$
then one can show  that $\bar{\nu}_{\gamma}\coloneqq \P_{\gamma}(\xi_t^{\{-\infty\}}\in\cdot)$ is an invariant law for the \textit{oriented percolation process} with parameter $\gamma$. We call $\bar{\nu}_{\gamma}$ the upper invariant measure for the process in the sense that, if $\nu$ is another invariant measure for the \textit{oriented percolation process}, then $\nu$ is stochastically dominated by $\bar{\nu}_{\gamma}$. Moreover, $\P_{\gamma}(\xi_{t}^{\Z}\in\cdot)$ converges weakly to $\bar{\nu}_{\gamma}$. We refer the reader to \cite{L} for more details.

We now write the statement of Theorem 1.1 in \cite{LS} for oriented percolation. In the following, $\nu_{\rho}$ denotes product measure on $\{0,1\}^{\Z^2}$ with density $\rho$.
\begin{theorem}[Ligget T.M. and Steif J.F., 2006]\label{Liggett} For all $\rho<1$, there exists $\gamma$ such that $\bar{\nu}_{\gamma}$ stochastically dominates $\nu_{\rho}$. In particular, $\P_{\gamma}(\xi_t^{\Z}\in\cdot)$ stochastically dominates $\nu_{\rho}$.
\end{theorem} 

We are now able to write the proofs of Lemma \ref{orien_perc} and Lemma \ref{compar}.

\vspace{0.5cm}

\hspace{-0.6cm}\textit{Proof of Lemma \ref{words}.} Define the sets
$$T_{m,1}=\{(u,v)\in\Z^2_+:3m\leq v < 4m;~ 4m\leq u+v < 16m \},$$
$$T_{m,2}=\{(u,v)\in\Z^2_+:3m\leq u < 4m;~ 4m\leq u+v < 16m\},$$
and consider the line segment passing by the points in $\mathcal{L}_m=\{(v_1,v_2)\in \Z^2:v_2=-v_1+m - 1\}$ (note that $L_{m}\subset\mathcal{L}_m$). For a mere technical reason, we will assume that the oriented percolation process takes place in $\LL^2$ instead of $\LL^2_+$ (this will be useful to obtain Equation \eqref{Lig} below).  With this in mind, let $$Q_{\mathcal{L}_{m}} = \{ x \in L_{16m} : \mathcal{L}_{m} \rightarrow x \text { in } \mathbb{L}^2 \}.$$

\begin{figure}[ht]
\centering
\begin{tikzpicture}[scale=0.03]

\draw[line width=0.03cm,black] (0,0) -- (0,200);
\draw[line width=0.03cm,black] (0,0) -- (200,0);
\draw[line width=0.03cm,black] (0,200) -- (50,150);
\draw[line width=0.03cm,black] (150,50) -- (200,0);
\draw[line width=0.09cm,black] (50,150) -- (150,50);

\draw[line width=0.09cm,black] (12.5,37.5) -- (37.5,12.5);
\draw[line width=0.03cm,gray] plot coordinates {(12.5,37.5)  (-20,70)};
\draw[line width=0.03cm,gray] plot coordinates {(37.5,12.5)  (70,-20)};
\draw[line width=0.03cm,dashed] plot coordinates {(-22,72)  (-32,82)};
\draw[line width=0.03cm,dashed] plot coordinates {(70,-20)  (82,-32)};

\draw[line width=0.05cm,lightgray, fill = lightgray] (37.5,12.5) -- (37.5,162.5) --  (50,150) --(50,0) -- (37.5,12.5);
\draw[line width=0.05cm,lightgray, fill = lightgray] (0,50) -- (12.5,37.5) -- (162.5,37.5) -- (150,50) --(0,50);

\draw[line width=0.05cm,dashed] plot coordinates {(30,20)  (70,20) (70,22) (130,22) (130,30) (170, 30)};
\draw[line width=0.03cm,black] plot  coordinates {(39,11) (39,47)(41,47) (41,100)(47,100)  (47,153)};
\draw[line width=0.03cm,black] plot  coordinates {(11,39) (48,39) (48,42)  (87,42) (87,48) (152,48)};

\draw (0,0)node[below] {$(0,0)$};
\draw (200,0)node[below] {$(16m,0)$};
\draw (105,57)node[right, gray] {$T_{m,1}$};
\draw (14,130)node[gray, right] {$T_{m,2}$};
\draw (115,120)node[below] {\normalsize $\mathbf{L_{16m,2}}$};
\draw (-35,57)node[right, black] {\normalsize $\mathcal{L}_{m}$};
\end{tikzpicture}
\caption{Graphical representation of the event $E_1\cap E_2\cap E_3$\label{event}.}
\end{figure}
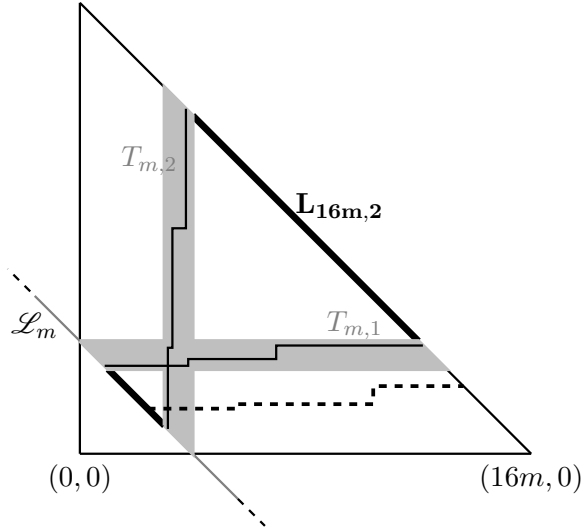

For $S\subset L_{4m,2}$ with $|S|\geq m$, consider the following events:
\begin{itemize}
\item $\mathrm{E_1}=\{\exists\, (u,v)\in L_{16m}\mbox{ such that } L_{4m,1}\rightarrow (u,v)\mbox{ in $T_{m,1}$}\}$,
\item $\mathrm{E_2}=\{S\rightarrow L_{16m}\}$,
\item $\mathrm{E_3}=\{\exists\, (u,v)\in L_{16m}\mbox{ such that } L_{4m,3}\rightarrow (u,v)\mbox{ in $T_{m,2}$}\}$,
\item $\mathrm{E}_4$= \{$Q_{\mathcal{L}_{4m}}\geq 4m$\}.
\end{itemize}
Observe that if $E_1\cap E_2\cap E_3$ occurs, then $M_S=Q_{\mathcal{L}_{4m}}$ (see Figure \ref{event}). If in addition $E_4$ occurs, then $M_S\geq 4m$, yielding
$$\P_{\gamma}\left(M_S< 4m\right)\leq \sum_{i=1}^4 \P_{\gamma}(E_i^c).$$

Let us bound the probabilities in the sum above. An upper bound for $\P_{\gamma}(E_1^c)$ can be obtained by a standard counting argument, which is simpler when considering oriented bond percolation. We will employ such argument to obtain an upper bound for $P_{q}(E'^c_{1})$, where ${E'_1}=\{\exists\, (u,v)\in L_{16m} \;:\; L_{4m,1}\rightarrow (u,v)\mbox{ in $T_{m,1}$}\}$ is the equivalent of $E_1$ for oriented bond percolation and $P_q(\cdot)$ denotes the product law of parameter $q$ for oriented bond percolation. We then use the fact $\P_{\gamma}(E_1) \geq P_{q}(E'_1)$ holds whenever $q = q(\gamma)$ is such that $1-\gamma = (1-q)^2$.

Pick a configuration $\omega \notin E'_1$, and note that the oriented connected component of $L_{m,1}$ in $T_{m,1}$ does not intersect $L_{16m,1}$ in $\omega$. This means there is a nearest-neighbor contour in the dual lattice, such that all sites immediately to the left of $\tau(\omega)$ belong to the open cluster of $L_{m,1}$, and no open edges starting at these sites cross $\tau(\omega)$. We also note that the \textit{west} and \textit{north} steps of this contour cannot cross open edges. For a fixed contour $\tau = (\tau_0, \dots, \tau_n)$, define $q_{\tau} = P_q(\tau(\omega) = \tau)$, and note that $q_{\tau} = (1-q)^{\alpha(\tau)}$, where $\alpha(\tau) = |\{1 \leq i \leq n : \tau_i - \tau_{i-1} \in \{(-1,0), (0,1)\}\}|$. Write $\tau_i = (x_i, y_i) \in \Z^2 +(1/2, 1/2)$ for $i = 0,\dots, n$, set $j_0 = 0$, and recursively define $j_{k+1} = \inf \{ i \in (j_k, n] : x_i > x_{j_k} \}$ for each $k$ such that the previous set is not empty. If $k$ is such that  $\{ i \in (j_k, n] : x_i > x_{j_k} \} = \emptyset$, we stop the recursion. 

Assume $\tau$ is given with associated sequence $(j_1, \dots, j_{\ell})$. Let $\sigma(\tau) = (\sigma_0, \dots, \sigma_n)$, with $\sigma_i = \tau_i - a_i(1,0)$, where $a_i$ is the cardinality of the set $\{ k \in [1,\ell] : j_k \leq i\}$. Note that if $\tau \neq \tau'$, then $\sigma(\tau) \neq \sigma(\tau')$. Writing $\sigma_i = (s_i, t_i)$ yields $s_i \leq s_0$, $i = 0, \dots, n$, which means $\sigma(\tau)$ is never on the right of $\sigma_0$. Moreover, a necessary condition for $\sigma_i = \sigma_{i+1}$ is that $s_i = s_0$.

Now let $\rho(\tau)$ denote the path obtained by deleting the repeated vertices of $\sigma(\tau)$. More precisely, consider the set $\{0\} \cup \{i \in [1,n] : \sigma_i \neq \sigma_{i-1}\}$, let $I(\tau) = (i_0, \dots, i_k)$ be the natural ordering of its elements, and set $\rho(\tau) = ( \sigma_{i_0}, \dots, \sigma_{i_k})$. Observe that $q_{\rho(\tau)} = q_{\tau}$, since in the construction of $\rho(\tau)$ only edges with \textit{east} direction are deleted. Since $\rho(\tau)$ is never on the right or bellow the first vertex $\sigma_{i_0}$, at least half of its edges are in the \textit{north} or \textit{west} directions, and thus $q_{\rho(\tau)} \leq (1 - q)^{k/2}$, where $k = |\rho(\tau)|$ is the size of $\rho(\tau)$. For each $j \in [0,k)$, let $z_j = i_{j+1} - i_j - 1$ denote the number of repetitions of $\sigma_{i_j}$ in $\sigma$. As observed before, if $z_j > 0$ then $s_{i_j} = s_0$. Also note that due to the size of the box $T_{m,1}$, we have $z_0 + \cdots + z_k \leq 12m$.

We say that $\rho$ is admissible if there exists $\omega \notin E'_1$ such that $\rho(\tau(\omega)) = \rho$. Denote by $\Lambda$ the set of admissible paths. Fix $\rho = (\rho_0, \dots, \rho_k) \in \Lambda$ and write $\rho_i = (x_i, y_i)$, $i \in [0,k]$. Also, write $R = \{ i \in [1,k] : x_i = x_0  \}$, and note that $|R|\leq m$. Given $\tau \in \Delta_{\rho}=\{\tau : \rho(\tau) = \rho \}$, note that the indices of $I(\tau) = (i_0, \dots, i_k)$ are such that $\rho_j = \sigma_{i_j}(\tau)$, for $i = 0, \dots, k$. According to these observations, it follows that $j\in R$ if $z_j(\tau) > 0$. Therefore, we have established a correspondence from $\Delta_{\rho}$ to the set $
\{(z_j)_{j \in R} : z_j \in \Z^+ \text{ and } \sum_{j \in R} z_j \leq 12m\}$. Since this correspondence is reversible, in the sense that is possible to obtain $\tau$ by adding edges with \textit{east} direction according to the sequence $(z_j)_{j \in R}$, such correspondence is injective. Therefore, for every $\rho \in \Lambda$, we have $|\Delta_{\rho}|\leq c_1^m$, for some universal constant $c_1$, independent of $m$ and $q$.

Finally, we obtain 
\begin{align*}
    \P_{\gamma}(E_1^c) &\leq \sum_{\rho \in \Lambda} \sum_{\tau \in \Delta_{\rho}}  (1 - q(\gamma))^{\alpha(\tau)} 
    \leq \sum_{\rho \in \Lambda} \sum_{\tau \in \Delta_{\rho}} (1 - q(\gamma))^{|\rho| / 2}
    \leq c_1^m \sum_{k \geq m}\sum_{\substack{\rho \in \Lambda \\ |\rho| = k}} (1 - q(\gamma))^{k/2} \\
    & \leq [c(1-q(\gamma))]^{m/2}, \quad \forall m \in \N,
\end{align*}
for some universal constant $c$, independent of $\gamma$. By symmetry the same bound holds for $\P_{\gamma}(E_3^c)$.

Consider now the event $E_2$. Note that, by symmetry, the probability of $\{S\rightarrow L_{16m}\}$ is the same as the probability of $\{L_{16m}\rightarrow S\}$ in the graph with reversed orientation, that is, the graph $(\Z^2,\E^*)$, where $\E^*=\{\langle v_1,v_2\rangle\subset \Z^2\times \Z^2: v_1\geq 0, v_2=v_1-(1,0)\mbox{ or } v_2=v_1-(0,1)\}$. An application of Theorem \ref{Liggett} shows that, given any $\beta\in(0,1)$, there exists $\gamma(\beta)\in(0,1)$ such that, for all $\gamma>\gamma(\beta)$,
$$\P_{\gamma}\left(L_{16m}\nrightarrow S\right)\leq P(Y=0),$$
where $Y\sim Bin(|S|,1-\beta)$. Therefore,
\begin{equation*}
\P_{\gamma}(E_2^c)\leq \beta^{|S|}\leq \beta^{m},\,\,\,\,\,\,\, \mbox{for all } \gamma>\gamma(\beta).
\end{equation*}

Take $\beta$ and $\gamma>\gamma(\beta)$ as above. A second application of Theorem \ref{Liggett} gives
\begin{equation}\label{Lig}
\P_{\gamma}(E_4^c)\leq P(X<4m),
\end{equation}
where $X\sim Bin(|L_{16m,2}|,1-\beta)=Bin(8m,1-\beta)$. Writing $X =8m-Z\sim Bin(8m,\beta)$ and applying a Chernoff bound, we obtain 
\begin{equation*}
\P(X<4m)=\P(Z>4m)\leq \left[\frac{(1-\beta+\beta e^t)^8}{e^{4t}}\right]^{m},\,\,\forall\, t\in\mathbb{R}.
\end{equation*}

Putting all together, for any $\beta\in (0,1)$, there exists $\gamma(\beta)\in(0,1)$ such that
$$\P_{\gamma}(M_S<4m)\leq 2[c(1-q(\gamma))]^{m/2}+\beta^{m}+\left[\frac{(1-\beta+\beta e^t)^8}{e^{4t}}\right]^{m},$$ for some universal constant $c>0$, for all $m\in \N$ and $\gamma>\gamma(\beta)$. Finally, given $a>0$, choose $t(a)$ sufficiently large, $\beta(t,a)$ sufficiently small, and  $\gamma(\gamma(\beta),c,a)$ such that
$$2[c(1-q(\gamma))]^{m/2}+\beta^{m}+\left[\frac{(1-\beta+\beta e^t)^8}{e^{4t}}\right]^{m}\leq a^{m},$$ and the result follows.
\qed

We now prove Lemma \ref{compar}.

\hspace{-0.6cm}\textit{Proof of Lemma \ref{compar}.} Fix $\eta\in\Xi_{16m}$. From the coupling described at the beginning of Section \ref{main_proof} and since \eqref{diverge} holds, it is clear that given any $p\in(0,1)$, $\epsilon>0$ and $\gamma\in(0,1)$, one can choose $N$ and $M$ sufficiently large so that the corresponding process of black points stochastically dominates an oriented independent percolation process on $\Z^2_+$ with density $\gamma\in(0,1)$ (see Equation \eqref{black} and Lemma 1 in \cite{GM}). This gives the desired result.
\qed

\section{Acknowledgements} The research of Pablo Gomes was partially supported by FAPESP, grant 2020/02636-3. Ot\'avio Lima was partially supported by CAPES. Roger Silva was partially supported by FAPEMIG, grants APQ-00774-21 and APQ-0868-21. This study was financed in part by the Coordena\c{c}\~{a}o de Aperfei\c{c}oamento de Pessoal de N\'{i}vel Superior, Brasil (CAPES), Finance Code 001. We also thank the author of the referee report.




\begin{thebibliography}{99}

\bibitem{BK} Benjamini I. and Kesten H.: Percolation of arbitrary words in $\{0,1\}^{\mathbb{N}}$, \emph{Ann. of Probab.}, {\bf 23}, (1995), 1024--1060.

\bibitem{CR} Campanino M. and Russo L.: An upper bound for the critical probability for the three-dimensional cubic lattice, \emph{Ann. of Probab.}, {\bf 13}, (1985), 478--491.

\bibitem{CCN} Chayes J.T., Chayes L. and Newman C.M.: Bernoulli percolation above threshold: and invasion percolation analysis, \emph{Ann. of Probab.}, {\bf 15}, (1987), 1272--1287.

\bibitem{D} Dekking F.M.: On the probability of occurrence of labelled subtrees of a randomly labelled tree, \emph{Theoret. Comput. Sci.}, {\bf 65}, (1989), 149--152.

\bibitem{ELV} van Enter A.C.D., de Lima B.N.B. and Valesin D.: Truncated long-range percolation on oriented graphs, \emph{J. Stat. Phys.}, {\bf 164}, (2016), 166--173.

\bibitem{FL} Friedli S. and de Lima B.N.B.: On the truncation of systems with non-summable interactions, \emph{Journ. Stat. Phys.}, {\bf 122}, (2006), 1215--1236.

\bibitem{GM} Grimmett G. and Marstrand J.M: The supercritical phase of percolation is well behaved, \emph{
Proceedings: Mathematical and Physical Sciences},
{\bf 430}, (1990), 439--457.


\bibitem{GS} Grimmett G. and Stacey A.M.: Critical probabilities for site and bond percolation models., \emph{Ann. of Probab.}, {\bf 26}, (1998), 1788--1812.

\bibitem{KSZ} Kesten H., Sidoravicius V. and Zhang Y.: Almost all words are seen in critical site percolation on the triangular lattice, \emph{Electron. J. Probab.}, {\bf 3}, (1998), 1--75.

\bibitem{KSZ2} Kesten H., Sidoravicius V. and Zhang Y.: Percolation of arbitrary words on the Close-Packed Graph of 
$\Z^2$, \emph{Electron. J. Probab.}, {\bf 6}, (2001), 1--27.

\bibitem{L} Liggett T.M.: Interacting Particle Systems, (2005), Springer Berlin.


\bibitem{LS} Liggett T.M. and Steif J.F.: Stochastic domination: the contact process, Ising models and FKG measures, \emph{Ann. Inst. H. Poincaré Probab. Statist.}, {\bf 42}(2), (2006), 223--243.

\bibitem{LSS} de Lima B.N.B., Sanchis R. and Silva R.W.C.: Percolation of words on $\Z^d$ with long-range connections, \emph{J. Appl. Probab.}, {\bf 48}(4), (2011), 1152--1162.

\bibitem{NTT} Nolin P., Teixeira A. and Tassion V.: No exceptional words for Bernoulli percolation, arXiv:1911.04816.

\bibitem{W} Wierman J.C.: AB percolation: a brief survey. \emph{Combinatorics and Graph Theory}, {\bf 25}, (1989), 241-251.


\end{thebibliography}
\end{document}